\documentclass{article}
\usepackage{geometry,enumerate}
\usepackage{amsmath, amssymb,amsthm,mathtools}
\usepackage[all]{xy}
\usepackage{graphics}
\usepackage{comment}
\usepackage{cleveref}
\usepackage{multicol}\columnseprule=0.4pt
\usepackage{color}
\title{A remark of the number of quasi-hereditary structures}
\author{Yuichiro Goto}
\date{}
\newtheorem{df}{Definition}[section]

\newtheorem{thm}[df]{Theorem}
\newtheorem{prop}[df]{Proposition}
\newtheorem{lem}[df]{Lemma}

\newtheorem{rmk}[df]{Remark}

\crefname{df}{Definition}{Definitions}
\crefname{eg}{Example}{Examples}
\crefname{thm}{Theorem}{Theorems}
\crefname{prop}{Proposition}{Propositions}
\crefname{lem}{Lemma}{Lemmas}
\crefname{cor}{Corollary}{Corollaries}
\crefname{rmk}{Remark}{Remarks}

\newcommand{\qhstr}{{\rm qh.str}}

\DeclareMathOperator{\module}{mod}

\DeclareMathOperator{\End}{End}

\DeclareMathOperator{\op}{op}

\begin{document}

\maketitle \vspace{3mm}

\begin{abstract}
Dlab and Ringel showed that algebras being quasi-hereditary in all total orders for indices of primitive idempotents becomes hereditary. 
So, we are interested in for which orders a given quasi-hereditary algebra is again quasi-hereditary. 
As a matter of fact, we consider permutations of indices instead of total orders.
If the standard modules defined by two permutations coincide, we say that the permutations are equivalent. 
Moreover if the algebra with permuted indices is quasi-hereditary, then this equivalence class of the permutation is called a quasi-hereditary structure.
In this article, we give a method of counting the number of quasi-hereditary structures for certain algebras.
\end{abstract} 

\section{Introduction}
Quasi-hereditary algebras, introduced by Cline, Parshall and Scott, generalize hereditary algebras.
Moreover Dlab and Ringel showed in Theorem 1 of \cite{dlab1989quasi} that if an algebra is quasi-hereditary in all orders, it becomes hereditary, and vice versa.
From this point of view, we study quasi-hereditary structures on a given algebra.
Flores, Kimura and Rognerud gave a method of counting the number of quasi-hereditary structures on a path algebras of Dynkin types in \cite{flores2022combinatorics}.
In this paper, we will use permutations instead of total orders when considering quasi-hereditary structures.

The followings are the main results of this article.
\begin{thm}[\Cref{result 1tu,result st1}]
    For some non-negative integers $s,t,u \geq 0$, let $Q(s,t,u)$ be the quiver
    
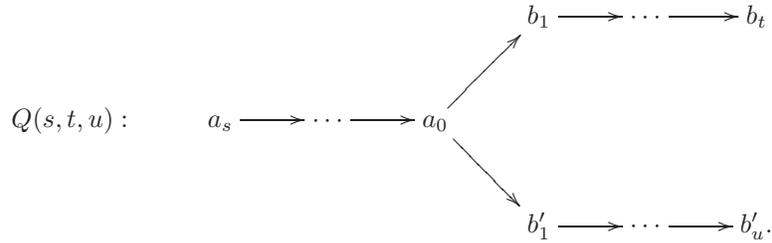
\[\xymatrix{
        &&&& b_1 \ar[r] & \cdots \ar[r] & b_t\\
        Q(s,t,u): & a_s \ar[r] & \cdots \ar[r] & a_0 \ar[ur]\ar[dr] \\
        &&&& b'_1 \ar[r] & \cdots \ar[r] & b'_u.
    }\]
    If one of $s,t,u$ is equal $1$, then the numbers of quasi-hereditary structures $q(s,t,u)$ on path algebras $KQ(s,t,u)$ are as follows.
    \begin{enumerate}[$(1)$]
        \item $q(1,t,u) = (c_{t+2} - c_{t+1})(c_{u+2} - c_{u+1})+c_{t+1}c_{u+1}$.
        \item $q(s,t,1) = q(s,1,t) = (t+1)c_{s+t+2} - 3tc_{s+t+1} - \sum_{k=0}^{t-2}(k+1)c_{t-k-1}c_{s+k+2}$
    \end{enumerate}
    Here, $c_k$ are the Catalan numbers.
\end{thm}


\section{Preliminaries}\label{Preliminaries}
Throughout this article, let $K$ be a field, $A$ a finite dimensional $K$-algebra with pairwise orthogonal primitive idempotents $e_1,\dots,e_n$, and let $\Lambda=\{1,\dots,n\}$.
For $i \in \Lambda$, we denote $P(i)=e_i A$ the indecomposable projective module, $S(i)$ the top of $P(i)$.
The category of finitely generated right $A$-modules is denoted by $\module A$ and call its objects just $A$-modules.
For non-negative integer $k \geq 0$, we write $c_k$ the Catalan number i.e. $c_k = \frac{(2k)!}{(k+1)! k!}$.

In this section, we recall the definition and the known properties of quasi-hereditary structures. 

\begin{df}\label{def of qha}
    Let $A$ be an algebra as above and $\sigma,\tau \in \mathfrak{S}_n$.
    \begin{enumerate}[$(1)$]
        \item For each $i \in \Lambda$, the $A$-module $\Delta^\sigma(i)$, called the {\bfseries standard module} with respect to $\sigma$, is defined by the maximal factor module of $P(i)$ having only composition factors $S(j)$ with $\sigma(j) \leq \sigma(i)$.
        Moreover we will write the set $\{\Delta^\sigma(1), \dots, \Delta^\sigma(n)\}$ by $\Delta^\sigma$. 
        \item We say that an $A$-module $M$ has a {\boldmath $\Delta^\sigma$}{\bfseries -filtration} 
        if there is a sequence of submodules 
        \[M=M_1 \supset M_2 \supset \dots \supset M_{m+1} = 0\] 
        such that for each $1 \leq k \leq m$, $M_k / M_{k+1} \cong \Delta^\sigma(j)$ 
        for some $j\in\Lambda$.
        \item A pair $(A,\sigma)$ is said to be a {\bfseries quasi-hereditary algebra} provided that the following conditions are satisfied.
        \begin{enumerate}[$(a)$]
            \item $\End_A(\Delta^\sigma(i))$ are local algebras for all $i \in \Lambda$.
            \item $A_A$ has a $\Delta^\sigma$-filtration.
        \end{enumerate}
        \item If $(A,\sigma)$ and $(A,\tau)$ are quasi-hereditary with $\Delta^{\sigma} = \Delta^{\tau}$, then the permutations $\sigma$ and $\tau$ are said to be {\bfseries equivalent}.
        Moreover, we call such a equivalence class by a {\bfseries quasi-hereditary structure} on $A$ and we denote by {\boldmath $\qhstr(A)$} the set of quasi-hereditary structures.
    \end{enumerate}
\end{df}

Now we recall the results in \cite{flores2022combinatorics}.

\begin{lem}[\cite{flores2022combinatorics} Lemma 2.22]\label{2.22}
    For the opposite algebra $A^{\op}$ of $A$, we have
    \[|\qhstr(A)| = |\qhstr(A^{\op})|.\]
\end{lem}

\begin{thm}[\cite{flores2022combinatorics} Corollary 4.8]\label{4.8}
    For a quiver $\xymatrix{Q_n:1 \ar[r] & 2 \ar[r] & \cdots \ar[r] & n}$ with $n \geq 1$, we have 
    \[|\qhstr(KQ_n)| = c_n.\]
\end{thm}

\begin{df}
    Let $Q$ be a finite quiver whose underlying graph is a tree. 
    A {\bfseries deconcatenation} of $Q$ at a sink or a source $v$ is a disjoint union $Q^1 \sqcup Q^2 \sqcup \cdots \sqcup Q^l$ of proper full subquivers $Q^i$ of $Q$ satisfying the following properties.
    \begin{enumerate}
        \item Each $Q^i$ is a connected full subquiver of $Q$ having a vertex $v$.
        \item On vertices, put $\overline{Q^i_0} = Q^i_0 \backslash \{v\}$.
        Then $Q_0 = \overline{Q^1_0} \sqcup \overline{Q^2_0} \sqcup \cdots \sqcup \overline{Q^l_0} \sqcup \{v\}$.
    \end{enumerate}
\end{df}

\begin{thm}[\cite{flores2022combinatorics} Theorem 3.7]\label{3.7}
    Let $Q$ be a finite quiver whose underlying graph is a tree and $Q^1 \sqcup Q^2 \sqcup \cdots \sqcup Q^l$ a deconcatenation of $Q$.
    Then we have 
    \[|\qhstr(KQ)| = \prod_{i=1}^{l} \left|\qhstr\left(KQ^i\right)\right|.\]
\end{thm}

\section{Main results}
By the results in the previous section, we next want to calculate the number of quasi-hereditary structures on the path algebra of a quiver
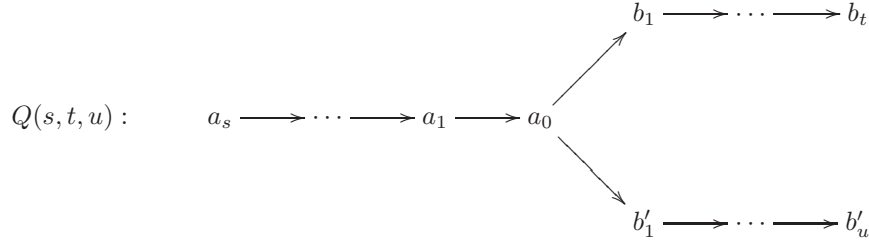
\[\xymatrix{
    &&&&& b_1 \ar[r] & \cdots \ar[r] & b_t\\
    Q(s,t,u): & a_s \ar[r] & \cdots \ar[r] & a_1 \ar[r] & a_0 \ar[ur]\ar[dr] \\
    &&&&& b'_1 \ar[r] & \cdots \ar[r] & b'_u
}\]
for some integers $s,t,u \geq 0$.
As a remark, there is an isomorphism $KQ(s,t,u) \cong KQ(s,u,t)$, in particular, $\qhstr(KQ(s,t,u)) = \qhstr(KQ(s,u,t))$.
So we may assume that $t \geq u$ without loss of generality.
Then Flores, Kimura and Rognerud showed the following formula.
\begin{prop}[\cite{flores2022combinatorics} Proposition 5.6]\label{5.6}
    \[\begin{split}
        |\qhstr(KQ(s,t,u))| =
        & \sum_{i=0}^{s}\left|\qhstr\left(\frac{KQ(s,t,u)}{\langle a_i \rangle}\right)\right| + \sum_{j=1}^{t}\left|\qhstr\left(\frac{KQ(s,t,u)}{\langle b_j \rangle}\right)\right| \\
        & + \sum_{k=1}^{u}\left|\qhstr\left(\frac{KQ(s,t,u)}{\langle b'_k \rangle}\right)\right| - \sum_{k=1}^{u}\sum_{j=1}^{t}\left|\qhstr\left(\frac{KQ(s,t,u)}{\langle b_i, b'_k \rangle}\right)\right|,
    \end{split}\]
    where $\langle e \rangle$ are two sided ideals $e KQ(s,t,u) e$ for any idempotents $e = a_i, b_j, b'_k$.
\end{prop}
Now we realize that there are isomorphisms 
\[\begin{split}
    \frac{KQ(s,t,u)}{\langle a_i \rangle} &\cong KQ(s-i,0,0) \times KQ(i-1,t,u)\\
    \frac{KQ(s,t,u)}{\langle b_j \rangle} &\cong KQ(0,t-j,0) \times KQ(s,j-1,u)\\
    \frac{KQ(s,t,u)}{\langle b'_k \rangle} &\cong KQ(0,0,u-k) \times KQ(s,t,k-1)\\
    \frac{KQ(s,t,u)}{\langle b_j,b'_k \rangle} &\cong KQ(0,t-j,0) \times KQ(0,0,u-k) \times KQ(s,j-1,k-1).
\end{split}\]
We write $q(s,t,u)$ to denote $|\qhstr(KQ(s,t,u))|$.
Then we get the following equality form \Cref{5.6}.
\begin{lem}
    \[\begin{split}
        q(s,t,u) 
        =& c_{s}c_{t}c_{u} + \sum_{i=1}^{s} c_{s-i} q(i-1,t,u) + \sum_{j=1}^{t} c_{t-j} q(s,j-1,u)\\
        & + \sum_{k=1}^{u} c_{u-k} q(s,t,k-1) - \sum_{j=1}^{t}\sum_{k=1}^{u} c_{t-j} c_{u-k} q(s,j-1,k-1).
    \end{split}\]
\end{lem}

By using this equality, in the case where the underlying graph of a quiver $Q$ is Dynkin type $\mathbb{A}$, $\mathbb{D}$ or $\mathbb{E}$, we can obtain the number $|\qhstr(KQ)|$ as follows.
\begin{lem}[\cite{flores2022combinatorics} Lemmas 5.8, 5.9, Example 5.10]
    The followings hold.
    \begin{enumerate}
        \item[$\mathbb{A}_n:$] $q(n-1,0,0) = c_{n}$ for $n \geq 1$ (\Cref{4.8}).
        \item[$\mathbb{D}_n:$] $q(n-3,1,1) = 2c_{n} - 3c_{n-1}$, $q(1,n-3,1) = 3c_{n-1} - c_{n-2}$ for $n \geq 3$.
        \item[$\mathbb{E}_6:$] $q(1,2,2) = 106$, $q(2,2,1) = 130$.
        \item[$\mathbb{E}_7:$] $q(1,3,2) = 322$, $q(2,3,1) = 416$, $q(3,2,1) = 453$.
        \item[$\mathbb{E}_8:$] $q(1,4,2) = 1020$, $q(2,4,1) = 1368$, $q(4,2,1) = 1584$.
    \end{enumerate}
\end{lem}

As generalizations of the above, we will describe $q(1,t,u)$ and $q(s,t,1)$ by using Catalan numbers.
\begin{thm}\label{result 1tu}
    We have the equality 
    \[q(1,t,u) = (c_{t+2} - c_{t+1})(c_{u+2} - c_{u+1})+c_{t+1}c_{u+1}.\]
    \begin{proof}
        We show that the equality holds by induction on $t+u$.
        For $t+u=0$, the equalities hold as $q(1,0,0) = c_{2}$.
        Assume that $t+u \geq 1$ and that $q(1,j,k) = (c_{j+2} - c_{j+1})(c_{k+2} - c_{k+1})+c_{j+1}c_{k+1}$ for $j+k < t+u$.
        Then we have the following equalities.

        \[\begin{split}
            q(1,t,u) 
            =& c_{1}c_{t}c_{u} + c_{1} q(0,t,u) + \sum_{j=1}^{t} c_{t-j} q(1,j-1,u) \\
            & + \sum_{k=1}^{u} c_{u-k} q(1,t,k-1) - \sum_{j=1}^{t}\sum_{k=1}^{u} c_{t-j} c_{u-k} q(1,j-1,k-1) \\ 
            =& c_{t}c_{u} + c_{t+1}c_{u+1} + \sum_{j=1}^{t} c_{t-j} ((c_{j+1} - c_{j})(c_{u+2} - c_{u+1})+c_{j}c_{u+1}) \\
            & + \sum_{k=1}^{u} c_{u-k} ((c_{t+2} - c_{t+1})(c_{k+1} - c_{k})+c_{t+1}c_{k}) \\
            & - \sum_{j=1}^{t}\sum_{k=1}^{u} c_{t-j} c_{u-k} ((c_{j+1} - c_{j})(c_{k+1} - c_{k})+c_{j}c_{k}) \\ 
            =& c_{t}c_{u} + c_{t+1}c_{u+1} + (c_{t+2}-2c_{t+1})(c_{u+2} - c_{u+1})+(c_{t+1}-c_{t})c_{u+1} \\
            & + (c_{t+2}-c_{t+1})(c_{u+2} - 2c_{u+1})+c_{t+1}(c_{u+1}-c_{u}) \\
            & - (c_{t+2}-2c_{t+1})(c_{u+2}-2c_{u+1}) + (c_{t+1}-c_{t})(c_{u+1}-c_{u}) \\ 
            =& c_{t+2}c_{u+2} - c_{t+2}c_{u+1} - c_{t+1}c_{u+2} + 2c_{t+1}c_{u+1} \\
            =& (c_{t+2} - c_{t+1})(c_{u+2} - c_{u+1})+c_{t+1}c_{u+1}.
        \end{split}\]
    \end{proof}
\end{thm}

In order to calculate $q(s,t,1)$, we use the following formulas.
Since (3) of the following probably does not appear in literatures, we give its proof here.
\begin{lem}\label{cat lem}
    We have the following equalities.
    \begin{enumerate}[$(1)$]
        \item $\frac{t+2}{2}c_{t+1} = (2t+1)c_{t}$.
        \item $\sum^{t}_{k=0} c_{k}c_{t-k} = c_{t+1}$.
        \item $\sum^{t}_{k=0} kc_{k}c_{t-k} = \frac{t}{2}c_{t+1}$.
    \end{enumerate}
    \begin{proof}
        $(1),(2)$ The equalities are well known.
        $(3)$ There is the equality 
        \[\sum^{t}_{k=0} kc_{k}c_{t-k} = \sum^{t}_{k=0} (t-k)c_{t-k}c_{k}.\]
        On the other hand, 
        \[\sum^{t}_{k=0} kc_{k}c_{t-k} + \sum^{t}_{k=0} (t-k)c_{t-k}c_{k} = \sum^{t}_{k=0} t c_{k}c_{t-k} = t c_{t+1}.\]
        Hence we conclude that $\sum^{t}_{k=0} kc_{k}c_{t-k} = \frac{t}{2}c_{t+1}$.
    \end{proof}
\end{lem}

\begin{thm}\label{result st1}
    We have the equality 
    \[
        q(s,t,1) = (t+1)c_{s+t+2} - 3tc_{s+t+1} - \sum^{t-2}_{k=0}(k+1)c_{t-k-1}c_{s+k+2}.
    \]
    \begin{proof}
        We show the statement by induction on $s+t$.
        For $s+t=0$, we have $q(0,0,1) = c_{2}$.
        Next, we assume that $s+t \geq 1$.
        For integers $j,k$ with $0 \leq k \leq j$, we define 
        \[\alpha^{j}_{k} =
        \begin{cases}
            j+1 &(k=j)\\
            -3j &(k=j-1)\\
            -(k+1)c_{j-k-1} &(0 \leq k \leq j-2).
        \end{cases}\]   
        Then we can write $q(i,j,1) = \sum^{j}_{k=0}\alpha^{j}_{k}c_{i+k+2}$ for $i+j < s+t$, by induction hypothesis.
        Hence there are the following equalities.{\small 
        \[\begin{split}    
            & q(s,t,1) \\
            =& c_{s}c_{t}c_{1} + \sum_{i=1}^{s} c_{s-i} q(i-1,t,1) + \sum_{j=1}^{t} c_{t-j} q(s,j-1,1) + c_{0} q(s,t,0) - \sum_{j=1}^{t} c_{t-j} c_{0} q(s,j-1,0)\\
            =& c_{s}c_{t} +\sum^{t}_{k=0} \alpha^{t}_{k} \sum_{i=1}^{s} c_{s-i}c_{i+k+1} + \sum_{j=1}^{t}\sum_{k=0}^{j-1} \alpha^{j-1}_{k} c_{t-j} c_{s+k+2} + c_{s+t+1} - \sum_{j=1}^{t} c_{t-j} c_{s+j}\\
            =& c_{t}c_{s} + \sum^{t}_{k=0} \alpha^{t}_{k} (c_{s+k+2} - \sum_{i=0}^{k+1} c_{s+i}c_{-i+k+1}) + \sum_{k=0}^{j-1}(\sum_{j=1}^{t} \alpha^{j-1}_{k} c_{t-j}) c_{s+k+2} + c_{s+t+1} - \sum_{j=1}^{t} c_{t-j} c_{s+j}\\
            =& c_{t}c_{s} + \sum^{t}_{k=0} \alpha^{t}_{k}c_{s+k+2} - \sum_{i=0}^{k+1} (\sum^{t}_{k=0}\alpha^{t}_{k}c_{-i+k+1}) c_{s+i} + \sum_{k=0}^{j-1} (\sum_{j=1}^{t}\alpha^{j-1}_{k} c_{t-j}) c_{s+k+2} + c_{s+t+1} - \sum_{j=1}^{t} c_{t-j} c_{s+j}.
        \end{split}\]}
        We focus on the coefficients of $c_{s+l}$ in the last expression.
        Obviously for $l < 0$ or $l > t+2$, the coefficients of $c_{s+l}$ are zero.
        For $0 \leq l \leq t+2$, the coefficients of $c_{s+l}$ can be obtained as follows.\\
        \begin{enumerate}
            \setlength{\leftskip}{8mm}
            \item[$c_{s+t+2}:$] $\alpha^{t}_{t} = t+1$. 
            \item[$c_{s+t+1}:$] $\alpha^{t}_{t-1} - \alpha^{t}_{t}c_{0} + \alpha^{t-1}_{t-1}c_{0} + 1
                = -3t$.
            \item[$c_{s+t}:$] $\alpha^{t}_{t-2} - \alpha^{t}_{t-1}c_0 - \alpha^{t}_{t}c_1 + \alpha^{t-1}_{t-2}c_0 + \alpha^{t-2}_{t-2}c_1 - c_{0}
                = - t + 1$.
            \item[$c_{s+l}:$] $\alpha^{t}_{l-2} - \sum_{k=l-1}^{t} \alpha^{t}_{k} c_{-l+k+1} + \sum_{j=l-1}^{t} \alpha^{j-1}_{l-2} c_{t-j} - c_{t-l}$, for $2 \leq l < t$.
            In order to look at this, we note the following.
            \[
                \sum_{k=l}^{t} \alpha^{t}_{k} c_{-l+k+1} = \sum_{j=l-1}^{t-1} \alpha^{j-1}_{l-2} c_{t-j}.
            \]
            Indeed,
            \[\begin{split}
                &\sum_{k=l}^{t} \alpha^{t}_{k} c_{-l+k+1} - \sum_{j=l-1}^{t-1} \alpha^{j-1}_{l-2} c_{t-j}\\
                =& (t+1)c_{t-l+1} - 3tc_{t-l} - \sum_{k=l}^{t-2}(k+1)c_{t-k-1}c_{k-l+1}\\
                &- (l-1)c_{t-l+1} + 3(l-1)c_{t-l} + \sum_{j=l+1}^{t-1}(l-1)c_{j-l}c_{t-j}\\
                =& (t-l+2)c_{t-l+1} - 3(t-l+1)c_{t-l} - \sum_{k=1}^{t-l-1}(k+l)c_{t-k-l}c_{k} + \sum_{j=1}^{t-l-1}(l-1)c_{j}c_{t-j-l}\\
                =& (t-l+2)c_{t-l+1} - 3(t-l+1)c_{t-l} - \sum_{k=1}^{t-l-1}(k+1)c_{t-k-l}c_{k}\\
                =& (t-l+2)c_{t-l+1} - 3(t-l+1)c_{t-l} - \left(\frac{t-l+2}{2}c_{t-l+1} - (t-l+2)c_{t-l}\right)\\
                =& \frac{t-l+2}{2}c_{t-l+1} - (2t-2l+1)c_{t-l}\\
                =& (2(t-l)+1)c_{t-l}- (2t-2l+1)c_{t-l}\\
                =& 0.
            \end{split}\]
            Thus we conclude that, for $2 \leq l < t$, the coefficient of $c_{s+l}$ equals to
            \[\begin{split}
                & \alpha^{t}_{l-2} - \sum_{k=l-1}^{t} \alpha^{t}_{k} c_{-l+k+1} + \sum_{j=l-1}^{t} \alpha^{j-1}_{l-2} c_{t-j} - c_{t-l}\\
                =& \alpha^{t}_{l-2} - (\alpha^{t}_{l-1} + \sum_{k=l}^{t} \alpha^{t}_{k} c_{-l+k+1}) + (\alpha^{t-1}_{l-2} + \sum_{j=l-1}^{t-1} \alpha^{j-1}_{l-2} c_{t-j}) - c_{t-l}\\
                =& \alpha^{t}_{l-2} - \alpha^{t}_{l-1} + \alpha^{t-1}_{l-2} - c_{t-l}\\
                =& -(l-1)c_{t-l+1}.
            \end{split}\]
            \item[$c_{s+1}:$] $- \sum_{k=0}^{t} \alpha^{t}_{k} - c_{t-1}$. 
            This equals to
            \[\begin{split}
                - \sum_{k=0}^{t} \alpha^{t}_{k} - c_{t-1}
                &= - (t+1)c_t + 3tc_{t-1} + \sum_{k=0}^{t-2}(k+1)c_{t-k-1}c_{k} - c_{t-1}\\
                &= - (t+1)c_t + 3tc_{t-1} + \frac{t+1}{2}c_{t} - tc_{t-1} - c_{t-1}\\
                &= - \frac{t+1}{2}c_{t} + 2tc_{t-1} - c_{t-1}\\
                &= - (2t-1)c_{t-1} + (2t-1)c_{t-1}\\
                &= 0.
            \end{split}\]
            \item[$c_{s}:$] $c_t - \sum_{k=0}^{t} \alpha^{t}_{k+1}$.
            This equals to
            \[\begin{split}
                c_t - \sum_{k=0}^{t} \alpha^{t}_{k+1}
                &= c_t - (t+1)c_{t+1} + 3tc_{t} + \sum_{k=0}^{t-2}(k+1)c_{t-k-1}c_{k+1}\\
                &= c_t - (t+1)c_{t+1} + 3tc_{t} + \frac{t}{2}c_{t+1} - tc_{t}\\
                &= - \frac{t+2}{2}c_{t+1} + (2t+1)c_{t}\\
                &= - (2t+1)c_{t} + (2t+1)c_{t}\\
                &= 0.
            \end{split}\]
        \end{enumerate}
        Summing up the above equalities, we obtain 
        \[q(s,t,1) = (t+1)c_{s+t+2} - 3tc_{s+t+1} - \sum^{t-2}_{k=0}(k+1)c_{t-k-1}c_{s+k+2},\] 
        as desired.
    \end{proof}
\end{thm}

\begin{rmk}
    Consider a quiver $Q$ whose underlying graph $\overline{Q}$ is the same as that of $Q(1,t,u)$ and $Q(s,t,1)$, i.e.
    \[\xymatrix{
        \overline{Q}:\cdot \ar@{-}[r] & \cdots \ar@{-}[r] & \cdot \ar@{-}[r] & \cdot \ar@{-}[r]\ar@{-}[d] & \cdot \ar@{-}[r] & \cdots \ar@{-}[r] & \cdot\\
        &&& \cdot
    }\]
    Then from \Cref{2.22,3.7} and our results, we can describe $|\qhstr(KQ)|$ by using Catalan numbers.
\end{rmk}


\begin{thebibliography}{9}
    \bibitem{dlab1989quasi} V. Dlab and C. M. Ringel. Quasi-hereditary algebras. {\it Illinois Journal of Mathematics}, 33(2):280--291, 1989.
    \bibitem{flores2022combinatorics} M. Flores, Y. Kimura, and B. Rognerud. Combinatorics of quasi-hereditary structures. {\it Journal of Combinatorial Theory, Series A}, 187:105559, 2022.
\end{thebibliography}

\end{document}